\documentclass[twoside]{article}
\usepackage{latexsym,amssymb,theorem,a4wide}
\usepackage{psboxit}
\usepackage{calc}

\usepackage[centertags,leqno]{amsmath}

\usepackage{enumerate}

\usepackage{epsfig}

\numberwithin{equation}{section}
\theoremstyle{plain}
\theorembodyfont{\itshape}
\newtheorem{lemma}{Lemma}[section]
\newtheorem{theorem}[lemma]{Theorem}

\newtheorem{proposition}[lemma]{Proposition}
\newtheorem{corollary}[lemma]{Corollary}
\newtheorem{definition}[lemma]{Definition}

\newtheorem{example}[lemma]{Example}
{\theorembodyfont{\rmfamily}
  \newtheorem{remark}[lemma]{Remark}
}

\theoremstyle{break}
{\theorembodyfont{\rmfamily}
  \newtheorem{remarkb}[lemma]{Remark}
}
\theoremstyle{plain}

\newcommand{\postskip}{\vspace{\theorempostskipamount}}

\newenvironment{proof}{\par\noindent\textit{\textbf{Proof: }}}%
{\hfill$\Box$\par\postskip}

\providecommand{\ml}[1]{}
\newcommand{\dup}{\mathrm d}

\usepackage{bbm}
\usepackage{amsmath}
\usepackage{ifthen}
\usepackage{calc}

\usepackage{color}
\definecolor{light}{gray}{.91}

\newlength{\SatzboxTextwidth}
\newlength{\GrauBoxTextwidth}
\newlength{\HauptsatzboxTextwidth}
\newlength{\hfillparboxwidth}
 \providecommand{\E}{{\ensuremath{\mathbf{E}}}}

 \providecommand{\N}{{\ensuremath{\mathbbm{N}}}}

 \renewcommand{\P}{{\ensuremath{\mathbf{P}}}}
 \providecommand{\P}{{\ensuremath{\mathbf{P}}}}
 
 \providecommand{\R}{{\ensuremath{\mathbbm{R}}}}

 \providecommand{\1}{{\ensuremath{\mathbbm{1}}}}

 \providecommand{\MCG}{{\ensuremath{\mathcal G}}}

\providecommand{\ru}[1]    {{{(#1)}}}
\providecommand{\rub}[1]   {{{\bigl(#1\bigr)}}}
\providecommand{\rubb}[1]  {{{\biggl(#1\biggr)}}}
\providecommand{\ruB}[1]   {{{\Bigl(#1\Bigr)}}}

\providecommand{\eckb}[1]  {{{\bigl[#1\bigr]}}}
\providecommand{\eckbb}[1] {{{\biggl[#1\biggr]}}}
\providecommand{\eckB}[1]  {{{\Bigl[#1\Bigr]}}}

\providecommand{\ddt}{{\ensuremath{\frac{d}{dt}}}}

\providecommand{\abs}[1]  {{\ensuremath{|#1|}}}
\providecommand{\absb}[1] {{\ensuremath{\bigl|#1\bigr|}}}

\providecommand{\absB}[1] {{\ensuremath{\Bigl|#1\Bigr|}}}

\providecommand{\al}      {{\ensuremath{\alpha}}}

\providecommand{\ld}      {{\ensuremath{\lambda}}}

\providecommand{\eps}     {{\ensuremath{\varepsilon}}}

\providecommand{\limN}{{\ensuremath{{\displaystyle \lim_{N \ra \infty}}}}}

\providecommand{\tlimsupN}  {\limsup_{N \ra \infty}}

\providecommand{\qasN} {\ensuremath{\quad\text{as }N\to\infty}}

\providecommand{\qqasN}{\ensuremath{\qquad\text{as }N\to\infty}}

\providecommand{\ra}{\rightarrow}

\providecommand{\lra}{\longrightarrow}

\providecommand{\wlim}{\Longrightarrow}

\providecommand{\wlimeps}%
      {{\ensuremath{\stackrel{\eps \rightarrow \infty}%
                            {\Longrightarrow}}}}

\providecommand{\mal}{\ensuremath{{\displaystyle \cdot}}}

\providecommand{\fa}{\ensuremath{\;\;\forall\;}}

\providecommand{\Gen}{{\ensuremath{\MCG}}}

\providecommand{\Law}[2][]{{\ensuremath{\mathcal L^{#1}\left(#2\right)}}}

\providecommand{\LawB}[2][]{{\ensuremath{\mathcal L^{#1}\Bigl(#2\Bigr)}}}

\providecommand{\ul}[1]{{\ensuremath{\underline{#1}}}}
\providecommand{\uline}[1]{{\ensuremath{\underline{#1}}}}

\providecommand{\clearemptydoublepage}%
    {\newpage{\pagestyle{empty}\cleardoublepage}}

\newlength{\mylen}

\newenvironment{beidschub}[1][5mm]{%
  \setlength{\mylen}{\textwidth}%
  \addtolength{\mylen}{-#1}%
  \addtolength{\mylen}{-#1}%
  \par\smallskip\noindent\hspace{#1}\begin{minipage}[t]{\mylen}}
  {\end{minipage}\par\smallskip\noindent}

\def\mathclapinternal#1#2{\clap{$\mathsurround=0pt#1{#2}$}}
\def\clap#1{\hbox to 0pt{\hss#1\hss}}

\def\mathclap{\mathpalette\mathclapinternal}

\makeatletter
\newlength{\Breit@}\newlength{\breit@}\newlength{\diff@renz}
\providecommand{\St@ckrel}[2]{%
  \settowidth{\Breit@}{\ensuremath{^{#1}}}
  \settowidth{\breit@}{\ensuremath{#2}}
  \ifthenelse{\Breit@>\breit@}{
    \setlength{\diff@renz}{(\Breit@-\breit@)/2}
    \hspace{\diff@renz}
  }{}}

\providecommand{\etStackrel}[2]{%
  \St@ckrel{#1}{#2}%
  &\stackrel{\mathclap{#1}}{#2}
  \St@ckrel{#1}{#2}}
\makeatother

\newcommand{\cadlag}{c\`adl\`ag}

\renewcommand{\P}{\mathbf{P}}
\begin{document}
%
%
%
%
%
%
%
\newcommand{\Title}{\uppercase{%
Graphical Representation Of Some Du\-ali\-ty Relations In Stochastic Population Mo\-dels
}}
\newcommand{\ShortTitle}{%
Graphical Representation Of Some Duality Relations
}
%
%
%
\newcounter{NumAuthor}
\setcounter{NumAuthor}{%
2 
}
\addtocounter{NumAuthor}{-1}
\newcommand{\AuthorOne}{\noindent\uppercase{%
Roland Alkemper
}}
\newcommand{\AuthorTwo}{\uppercase{%
Martin Hutzenthaler
}}
\newcommand{\AuthorThree}{\uppercase{%
}}
\newcommand{\AuthorFour}{\uppercase{%
}}
\newcommand{\AuthorFive}{\uppercase{%
}}
%
%
%
\newcommand{\AddressOne}{%
Department of Mathematics,
Johannes-Gutenberg Universit\"at,
Staudingerweg 9,
55099 Mainz,
Germany}
\newcommand{\AddressTwo}{%
Department of Mathematics and Computer Science,
Universit\"at Frankfurt,
Robert-Mayer-Str. 6-10,
60325 Frankfurt am Main,
Germany}
\newcommand{\AddressThree}{%
}
\newcommand{\AddressFour}{%
}
\newcommand{\AddressFive}{%
}
%
%
%
\newcommand{\EmailOne}{%
alkemper@mathematik.uni-mainz.de
}
\newcommand{\EmailTwo}{%
hutzenth@math.uni-frankfurt.de
}
\newcommand{\EmailThree}{%
email.com                               
}
\newcommand{\EmailFour}{%
email.com                               
}
\newcommand{\EmailFive}{%
email.com                               
}
%
%
%
%
%
\newcommand{\Submitted}{February 22, 2007}
\newcommand{\Accepted}{XXX}
%
%
%
\newcommand{\SubjectClassification}{%
60K35
}
\newcommand{\Keywords}{%
Duality, graphical representation, 
Feller's branching diffusion,
branching-coalescing particle process, resampling-selection model,
stochastic population dynamics
}
%
%
%
\newcommand{\Abstract}{%
  \noindent
  We derive a unified
  stochastic picture for the duality of a resampling-selection
  model with a branching-coalescing particle
  process (cf.~\cite{AS05})
  and for the self-duality of
  Feller's branching diffusion with logistic growth
  (cf.~\cite{HW07}).
  The two dual processes are approximated by
  particle processes which are forward and backward processes in a
  graphical representation.
  We identify duality relations between the basic building blocks of the
  particle processes which lead to the two dualities mentioned above.
}
%
%
%
%
\thispagestyle{plain}
\rule{0in}{0pt}

%
%
\vskip 60pt
\noindent
{\Large\bf\uppercase{\Title}}

%
%
\vskip 12pt
\ifcase\value{NumAuthor}
	\uppercase{\AuthorOne}\\
	{\em \AddressOne}\\
	Email:\hskip 4pt\texttt{\EmailOne}
\or
	\uppercase{\AuthorOne}\\
	{\em \AddressOne}\\
	Email:\hskip 4pt\texttt{\EmailOne}\\[2mm]
	\uppercase{\AuthorTwo}\\
	{\em \AddressTwo}\\
	Email:\hskip 4pt\texttt{\EmailTwo}
\or
	\uppercase{\AuthorOne}\\
	{\em \AddressOne}\\
	Email:\hskip 4pt\texttt{\EmailOne}\\[2mm]
	\uppercase{\AuthorTwo}\\
	{\em \AddressTwo}\\
	Email:\hskip 4pt\texttt{\EmailTwo}\\[2mm]
	\uppercase{\AuthorThree}\\
	{\em \AddressThree}\\
	Email:\hskip 4pt\texttt{\EmailThree}
\or
	\uppercase{\AuthorOne}\\
	{\em \AddressOne}\\
	Email:\hskip 4pt\texttt{\EmailOne}\\[2mm]
	\uppercase{\AuthorTwo}\\
	{\em \AddressTwo}\\
	Email:\hskip 4pt\texttt{\EmailTwo}\\[2mm]
	\uppercase{\AuthorThree}\\
	{\em \AddressThree}\\
	Email:\hskip 4pt\texttt{\EmailThree}\\[2mm]
	\uppercase{\AuthorFour}\\
	{\em \AddressFour}\\
	Email:\hskip 4pt\texttt{\EmailFour}
\or
	\uppercase{\AuthorOne}\\
	{\em \AddressOne}\\
	Email:\hskip 4pt\texttt{\EmailOne}\\[2mm]
	\uppercase{\AuthorTwo}\\
	{\em \AddressTwo}\\
	Email:\hskip 4pt\texttt{\EmailTwo}\\[2mm]
	\uppercase{\AuthorThree}\\
	{\em \AddressThree}\\
	Email:\hskip 4pt\texttt{\EmailThree}\\[2mm]
	\uppercase{\AuthorFour}\\
	{\em \AddressFour}\\
	Email:\hskip 4pt\texttt{\EmailFour}\\[2mm]
	\uppercase{\AuthorFive}\\
	{\em \AddressFive}\\
	Email:\hskip 4pt\texttt{\EmailFive}
\fi
%
%
\vskip 6pt

%
%

\vskip 6pt\noindent
AMS 2000 Subject classification: \SubjectClassification

%
%

\vskip 18pt\noindent
\emph{Abstract \vskip 2pt } \Abstract

%
%
\section{Introduction}\label{sec:introduction}
Two processes $(X_t)_{t\geq0}$ and $(Y_t)_{t\geq0}$
with state spaces $E_1$ and $E_2$, respectively, are called dual with
respect to the duality function $H$ if $H\colon E_1\times E_2\to\R$
is a measurable and bounded function and if
$\E^x[ H(X_t,y) ]=\E^y[ H(x,Y_t) ]$ holds for all $x\in E_1$, $y\in E_2$
and all $t\geq0$ (see e.g. \cite{Lig85}).
Here superscripts as in $\P^x$ or in $\E^x$ indicate the initial
value of a process.
In this paper, $E_1$ and $E_2$ will be subsets of $[0,\infty)$
or will be equal to $\{0,1\}^N$.
We speak of a \emph{moment duality}
if $H(x,y)=y^x$ or $H(x,y)=(1-y)^x$,
$x\in E_1 \subset \N_0$, $y\in [0,1]$, and of a
\emph{Laplace duality} if
$H(x,y)=\exp\ru{-\ld x\mal y}$, $x,y\in E_1=E_2\subset [0,\infty)$,
for some $\ld>0$.

We provide a unified stochastic picture for the following moment duality
and the following Laplace duality of prominent processes from the field
of stochastic population dynamics.
For the moment duality, let $b,c,d\geq0$.
Denote by $X_t\in\N_{0}$ the number of particles
at time $t\geq0$ of
the branching-coalescing particle process
defined by the initial value $X_0=n$ and
the following dynamics: Each particle splits
into two particles at rate $b$, each particle dies at
rate $d$ and each ordered pair of particles coalesces into
one particle at rate $c$. All these events occur independently of each other. In the
notation of Athreya and Swart \cite{AS05}, this is the
$(1,b,c,d)$-braco-process.
Its dual process $(Y_t)_{t\geq0}$ is the unique strong solution  with
values in $[0,1]$ of the one-dimensional stochastic
differential equation
\begin{equation}  \label{eq:resem}
  \dup Y_t=(b-d)Y_t\,\dup t-b Y_t^2\,\dup t+\sqrt{2c Y_t(1-Y_t)}\,\dup B_t,
  \quad Y_0=y,
\end{equation}
\enlargethispage{2mm}%
where $(B_t)_{t\geq0}$ is a standard Brownian motion.
Athreya and Swart~\cite{AS05}
call this process the re\-samp\-ling-se\-lec\-tion process with
selection rate $b$, resampling rate $c$ and mutation rate $d$,
or shortly the $(1,b,c,d)$-resem-process. They prove
the moment duality
\begin{equation} \label{eq:braco_resem}
  \E^n\eckb{(1-y)^{X_t}} =\E^y\eckb{(1-Y_t)^n}
  \qquad \forall \, n\in\N_0,\, y\in[0,1],\, t\geq0.
\end{equation}
For the Laplace duality, let
$(X_t)_{t\geq0}$
denote Feller's branching diffusion with logistic growth,
i.e., the strong solution of
\begin{equation}  \label{eq:logistic_Feller}
  \dup X_t= \al X_t\,\dup t- \gamma X_t^2\,\dup t+\sqrt{2\beta X_t}\,\dup B_t,
\end{equation}
where $\al,\gamma,\beta\geq 0$ and 
$(B_t)_{t\geq0}$ is a standard Brownian motion.
We call this process the logistic Feller diffusion with parameters
$(\al,\gamma,\beta)$. 
Let $(Y_t)_{t\geq0}$ be a logistic Feller diffusion
with parameters $(\al,r\beta,\gamma/r)$ for some $r>0$.
Hutzenthaler and Wakolbinger~\cite{HW07} establish the Laplace duality
\begin{equation}  \label{eq:Fellog_dual}
     \E^x \eckb{ e^{- r X_t\cdot y}}
    =\E^y \eckb{ e^{- r x\cdot Y_t}},
     \qquad  \forall\, x,y\in[0,\infty),\, t\geq0.
\end{equation}
The duality relations~\eqref{eq:braco_resem} and~\eqref{eq:Fellog_dual}
include as special cases (see Remark~\ref{r:GWdual}
and Remark~\ref{r:Felldual})
the Laplace duality of Feller's branching diffusion with a
deterministic process,
the moment duality of the Fisher-Wright diffusion
with Kingman's coalescent, and the moment duality of the (continuous time)
Galton-Watson process with a deterministic process.

In the references~\cite{AS05} and~\cite{HW07},
the duality relations~\eqref{eq:braco_resem} and~\eqref{eq:Fellog_dual}
are proved analytically by means of a generator calculation.
In this paper, we take a different approach by explaining the dynamics of
the processes via \emph{basic mechanisms} on
the level of particles which lead to the above dualities.
To this end, for every $N \in \N$, we construct approximating Markov processes 
$\rub{X_t^N}_{t\geq0}$
and $\rub{Y_t^N}_{t\geq0}$ with \cadlag\ sample paths and
state space $\{0,1\}^N$ and with the following properties. 
The processes $\ru{X_t^N}_{t\geq0}$ and
$\ru{Y_t^N}_{t\geq0}$ are dual in the sense that
\begin{equation}   \label{eq:wedge_duality}
  \P^{x^N}\eckb{ X_t^N\wedge  y^N=\ul{0}}
  =\P^{y^N}\eckb{ x^N\wedge Y_t^N=\ul{0}},
  \quad\fa x^N, y^N\in\{0,1\}^N\,\fa t\geq0.
\end{equation}
The notation $x^N\wedge y^N$ denotes
component-wise minimum and $\ul{0}$ denotes the zero configuration.
If $\abs{X_0^N}=n$, for some fixed $n\leq N$,
then $\rub{\abs{X_t^N}}_{t\geq0}$ converges weakly to a
branching-coalescing particle process as $N\to\infty$.
We use the notation $\abs{x^N}:=\sum_{i=1}^N x^N_i$ for $x^N\in\{0,1\}^N$.
Assume that the set
of \cadlag-paths is equipped with the Skorohod topology (see e.g. \cite{EthKu}).
If $n=n(N)$ depends on $N$ such that
$n/N\to x\in[0,1]$ as $N\to\infty$,
then $\ru{\abs{X_t^N}/N}_{t\geq0}$ converges weakly  to
a resampling-selection model. If $n=n(N)$ satisfies
$n/\sqrt{N}\to x\geq0$, then
$\rub{\abs{X_{t\sqrt{N}}^N}/\sqrt{N}}_{t\geq0}$
converges weakly to Feller's branching diffusion with logistic growth.
The process $\ru{Y_t^N}_{t\geq0}$ differs
from $\ru{X_t^N}_{t\geq0}$ only by the set of parameters and
by the initial condition.

We will derive the moment duality~\eqref{eq:braco_resem}
and the Laplace duality~\eqref{eq:Fellog_dual}
from~\eqref{eq:wedge_duality} in the following way.
Let the random variable $X_0^N$ be uniformly distributed over all
configurations $x^N\in\{0,1\}^N$ with total number of individuals
of type $1$ equal to
$\abs{x^N}=n=n(N)$ for a given $n(N)\leq N$. Similarly, choose
$Y_0^N$ uniformly in $\{0,1\}^N$ 
with $\abs{Y_0^N}=k=k(N)$ for a given $k(N)\leq N$.
We will prove in Proposition~\ref{basicduality} that
property~\eqref{eq:wedge_duality} implies a prototype duality relation, namely
\begin{equation}  \label{dadudi}
   \limN\E\eckB{1-\frac{k}{N}}^{\absb{X_{t T_N}^N}}
  =\limN\E\eckB{1-\frac{\absb{Y_{tT_N}^N}}{N}}^{n},\quad t\geq0,
\end{equation}
under some assumptions -- including the convergence of both sides --
on the two processes and
on the sequence $\ru{T_N}_{N\geq1}\subset\R_{\geq0}$.
Choosing $n$ fixed, $k$ such that $\tfrac{k}{N}\to y\geq0$ and letting $T_N=1$,
we deduce from~\eqref{dadudi} 
(and from the convergence properties of $(X_t^N)_{t\geq0}$ and
of $(Y_t^N)_{t\geq0}$) the moment duality of 
a branching-coalescing particle process
with 
a resampling-selection model (cf.\ Theorem~\ref{GWdual}).
In order to obtain a Laplace duality of logistic Feller diffusions,
choose $n,k$ such that
$\tfrac{n}{\sqrt{N}}\to x\geq0$,
$\tfrac{k}{\sqrt{N}}\to y\geq0$
and $T_N=\sqrt{N}$.
Notice that $\ru{1-\tfrac{y}{\sqrt{N}}}^{x\sqrt{N}}$
converges to $e^{-xy}$ uniformly in $0\leq x,y\leq\tilde{x}$ as $N\to\infty$
for every $\tilde{x}\geq0$.
This together with the weak convergence of the rescaled processes 
will imply
\begin{equation}
  \limN\E \eckB{e^{- \abs{X^N_{t \sqrt{N}}}\cdot y \big/ \sqrt{N}}}
 =\limN\E \eckB{e^{-x \cdot \abs{Y^N_{t \sqrt{N}}}\big/ \sqrt{N}}} .
\end{equation}

\begin{figure}[ht]
\begin{beidschub}[0.5cm]
\begin{minipage}[t]{\linewidth}
  \epsfig{file=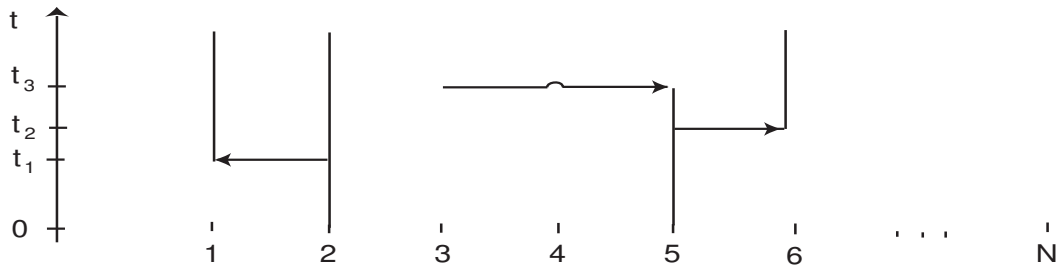, clip=,
      width=\linewidth, height=3.5cm}
  \caption{\footnotesize Three resampling events.
       Type $1$ is indicated by black lines,
       absent lines correspond to type 0.}
   \label{f:moran}
\end{minipage}
\end{beidschub}
\end{figure}

For the construction of the approximating processes,
we interpret the elements of $\{1,\ldots,N\}$ as ``individuals'' and
the elements of $\{0,1\}$ as the ``type'' of an individual. 
In the terminology of population genetics, individuals are
denoted as ``genes'', whereas
in population dynamics, the statement ``individual $i$ is of type $1$ (resp. $0$)''
would be
phrased as ``site $i$ is occupied (resp. not occupied) by a particle''. 
Throughout the paper, we assume that 
whenever a change of the configuration happens at most two individuals are involved.
We call every function $f\colon \{0,1\}^2\to\{0,1\}^2$ a
\textit{basic mechanism}.
A finite tuple $(f_1,...,f_m)$, $m\in\N$, of basic mechanisms together
with rates $\ld_1,...,\ld_m\in[0,\infty)$ defines a process with
state space $\{0,1\}^N$ by means of the following graphical representation,
which is in the spirit of Harris~\cite{Ha78}.
With every $k\leq m$ and 
every ordered pair $(i,j)\in\{1,...,N\}^2$, $i\neq j$, of individuals,
we associate a
Poisson process with rate parameter $\ld_k$.
At every time point of this Poisson process, the configuration
of $(i,j)$ changes according to $f_k$.
For example, if the pair of types was $(1,0)$ before,
then it changes to $f_k\ru{1,0}\in\{0,1\}^2$.
All Poisson processes are supposed to be independent.
This construction can be visualised by drawing arrows from $i$ to $j$ at
the time points of the Poisson processes associated with
the pair $(i,j)$ (cf.\ Figure~\ref{f:moran}).

As an example, consider the following continuous time
\emph{Moran model} $(M_t^N)_{t\geq0}$ with state space $\{0,1\}^N$.
This is a population genetic model where ordered pairs of individuals resample at rate $\beta/N$, $\beta>0$.
When a resampling event occurs at $(i,j)$,
individual $i$ bequeaths
its type to individual $j$.
Thus, the basic mechanism is $f^R$ defined by
\begin{equation}  \label{eq:fR}
  f^R(1,\cdot):=(1,1),\ f^R(0,\cdot):=(0,0).
\end{equation}
Figure~\ref{f:moran} shows a realisation with three resampling events.
At time $t_1$, the pair $(2,1)$ resamples.
The arrow in Figure~\ref{f:moran} at time $t_1$ indicates
that individual $2$ bequeaths its type to individual $1$.
Furthermore, individual $5$ inherits the type
of individual $3$ at time $t_3$.
The dual process of the Moran model is a coalescent process.
This process is defined by the coalescent mechanism $f^C$ given by
\begin{equation}  \label{eq:fC}
  f^C(1,\cdot):=(0,1),\ f^C(\uline{z}):=\uline{z},\quad \uline{z}\in\{(0,0),(0,1)\},
\end{equation}
and  by the rate $\beta/N$.
To put it differently,
the coalescent process is a coalescing
random walk on the complete oriented graph of $\{1,\ldots,N\}$.
In Section~\ref{sec:dual_basic_mechanisms}, we will specify
in which sense $f^R$ and $f^C$ are dual, and why this
implies~\eqref{eq:wedge_duality} (see Proposition~\ref{dualcrit}).
More generally, we will identify all dual pairs of basic mechanisms.

Our method elucidates the role of the square
in~\eqref{eq:logistic_Feller} for the duality of the
logistic Feller diffusion with another logistic Feller diffusion.
We illustrate this by the Laplace duality of
Feller's branching
diffusion $(F_t)_{t\geq0}$, which is the
logistic Feller diffusion with parameters $(0,0,\beta)$, $\beta>0$.
Its dual process $(y_t)_{t\geq0}$ is
the logistic Feller diffusion with parameters $(0,\beta,0)$,
i.e., the solution of the ordinary differential equation
\begin{equation}   \label{eq:def:y_t}
  \ddt y_t=-\beta\, y_t^2,\qquad y_0=y \in [0,\infty).
\end{equation}
The duality relation between these two processes is 
  $\E^x [e^{-F_ty}]=e^{-xy_t}$, $t \geq 0$.
In Theorem~\ref{Felldual}, we prove that the rescaled Moran model
$\rub{\abs{M_{t\sqrt{N}}^N}/\sqrt{N}}_{t\geq0}$
converges weakly to $\ru{F_t}_{t\geq0}$ as $N\to\infty$.
To get an intuition for this convergence, notice that
$(\abs{M_t^N})_{t\geq0}$ is a pure birth-death process
with size-dependent transition rates
(``birth'' corresponds to creation of an individual with type $1$,
whereas ``death'' corresponds to creation of an individual with
type $0$).
It remains to prove that the birth and death events
become asymptotically independent as $N\to\infty$.
It is known, see e.g.~Section 2 in \cite{Don84},
that the dual process of the Moran model $(M_t^N)_{t\geq0}$, $N\geq1$,
is a coalescing random walk. Furthermore,
the total number of particles of
this coalescing random walk is a pure death process
on $\{1,...,N\}$ which jumps from $k$
to $k-1$ at exponential rate $\tfrac{\beta}{N}k(k-1)$, $2\leq k\leq N$.
This rate is essentially quadratic in $k$ for large $k$.
We will see
that a suitably rescaled pure death process
converges to a solution of~\eqref{eq:def:y_t};
see Remark~\ref{r:pure_death_process}.
The square in~\eqref{eq:def:y_t}
originates in
the quadratic rate of the involved pure death process;
see the equations~\eqref{eq:puma} and~\eqref{eq:transition_rates} for details.

In the literature, e.g.\ \cite{Lig85},
the duality function $H(x^N,y^N)=\1_{x^N\leq y^N}$, $x^N,y^N\in\{0,1\}^N$,
can be found frequently, where
$x^N\leq y^N$ denotes component-wise comparison.
Processes $(X_t^N)_{t\geq0}$ and $(Y_t^N)_{t\geq0}$ with state space
$\{0,1\}^N$ are dual with respect to this duality
function if they satisfy
\begin{equation}   \label{eq:subset_duality}
  \P^{x^N}\eckb{X_t^N\leq y^N}=
  \P^{y^N}\eckb{x^N\leq Y_t^N}
  \quad\fa x^N,y^N\in\{0,1\}^N,t\geq0.
\end{equation}
The biased voter model is dual to a coalescing
branching random walk in this sense
(see~\cite{KN97}).
Property~\eqref{eq:subset_duality} could also be used to
derive the dualities mentioned in this introduction.
In fact, the two properties~\eqref{eq:wedge_duality}
and~\eqref{eq:subset_duality} are equivalent in the following
sense: If $(X^N_t)_{t\geq0}$ and $(Y^N_t)_{t\geq0}$
satisfy~\eqref{eq:wedge_duality}
then $(X^N_t)_{t\geq0}$  and $(\uline{1}- Y^N_t)_{t\geq0}$
satisfy~\eqref{eq:subset_duality} and vice versa.
In the configuration $\uline{1}$ every individual has type 1
and $\uline{1}-y$ denotes component-wise subtraction.
The dynamics of the process
$\ru{\uline{1}- Y_t^N}_{t\geq0}$ is easily
obtained from the dynamics of $(Y^N_t)_{t\geq0}$ by interchanging
the roles of the types $0$ and $1$.

%
%
\section{Dual basic mechanisms}
\label{sec:dual_basic_mechanisms}%
Fix $m\in\N$ and let $(X_t^N)_{t\geq0}$ and $(Y_t^N)_{t\geq0}$
be two processes defined by basic mechanisms
$(f_1,...,f_m)$ and $(g_1,...,g_m)$, respectively.
Suppose that the Poisson processes associated with $k\leq m$ have the
same rate parameter $\ld_k\geq0$, $k=1,\ldots,m$.
We introduce a property of basic mechanisms which will
imply~\eqref{eq:wedge_duality}.

%
%
\begin{definition}\label{dualfunc}
Let $f,g:\{0,1\}^2 \rightarrow \{0,1\}^2$ and for
$x=(x_1,x_2)\in \{0,1\}^2$ let $x^\dagger := (x_2,x_1)$.
The basic mechanisms $f$ and $g$ are said to be \textbf{dual} iff the following two conditions hold:
\begin{eqnarray}
  \forall \, x,y \in \{0,1\}^2\colon \;\;
       y \wedge \rub{f(x)}^\dagger = (0,0) \,
  &\Longrightarrow&g(y) \wedge x^\dagger = (0,0), \label{eq:f_bestimmt_g}\\
  \forall \, x,y \in \{0,1\}^2\colon\; \;
       x \wedge \rub{g(y)}^\dagger = (0,0) \,
  &\Longrightarrow& f(x) \wedge y^\dagger = (0,0).\label{eq:g_bestimmt_f}
\end{eqnarray}
\end{definition}
To see how this connects to the
duality relation in~\eqref{eq:wedge_duality}, we illustrate 
this definition by an example.
\begin{example}  \label{ex:1}
  \upshape%
  The resampling mechanism $f^R$ defined in~\eqref{eq:fR} and
  the coalescent mechanism $f^C$ defined in~\eqref{eq:fC} are dual.
  We check condition~\eqref{eq:f_bestimmt_g}
  with $f=f^{R}$ and $g=f^C$
  by looking at Figure~\ref{f:bd_bestimmt_c}.
  \newlength{\dummy}%
  \setlength{\dummy}{\linewidth*\real{0.8}}%
  \newlength{\reinschub}%
  \setlength{\reinschub}{(\linewidth-\dummy)/2}%
  \begin{figure}[ht]
  \begin{beidschub}[\reinschub]
  \noindent
  \begin{minipage}[t]{\linewidth}
    \epsfig{file=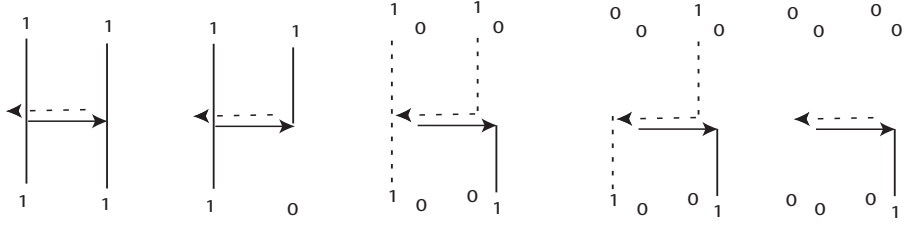, clip=,
          width=\dummy, height=3.0cm}
    \caption{\footnotesize The resampling mechanism 
       and the coalescent mechanism
       satisfy~\eqref{eq:f_bestimmt_g}}%
     \label{f:bd_bestimmt_c}
  \end{minipage}
  \end{beidschub}
  \end{figure}
  \noindent
  The resampling mechanism acts
  in upward time (solid lines),
  the coalescent mechanism in downward time (dashed lines).
  There are three nontrivial configurations for $x$, i.e., $(1,1)$, $(1,0)$
  and $(0,1)$. In the first two cases, we have $f^{R}(x)=(1,1)$.
  Then only $y=(0,0)$ satisfies $y\wedge \ru{f^{R}(x)}^\dagger=(0,0)$.
  In the third case, every $y$ satisfies
  $y\wedge \ru{f^{R}(0,1)}^\dagger=(0,0)$ and has to be checked separately.
  We see that whenever the configuration $y$ is
  disjoint from $\ru{f(x)}^\dagger$, i.e.,
  $y\wedge \ru{f(x)}^\dagger=(0,0)$,
  then $g(y)$ is disjoint from $x^\dagger$.
  The coalescent mechanism is the natural dual mechanism of
  the resampling mechanism.
  Type $1$ of the coalescent mechanism
  ``traces back''
  the lines of descent of type $0$ of the resampling mechanism.
  The ``birth event'' $(0,1)\mapsto(0,0)$ of an individual of type $0$
  results in a coalescent event of ancestral lines.

  Figure~\ref{f:c_bestimmt_bd} is useful to verify
  condition~\eqref{eq:g_bestimmt_f}.
  Again, the coalescent mechanism is drawn with dashed lines.
  Here, the coalescent process is started in the nontrivial
  configurations
  $(1,1)$, $(1,0)$ and $(0,1)$.
  In any case we obtain $\ru{f^{C}(y)}^\dagger=(1,0)$.
  Hence, all admissible $x$ are of the form $(0,\cdot)$.
  Condition~\eqref{eq:g_bestimmt_f} then follows from 
  $f^{R}(0,\cdot)=(0,0)$.
  \newlength{\duneu}%
  \setlength{\duneu}{\linewidth*\real{0.8}}%
  \newlength{\neuschub}%
  \setlength{\neuschub}{(\linewidth-\duneu)/2}%
  \begin{figure}[ht]
  \begin{beidschub}[\neuschub]
  \noindent
  \begin{minipage}[t]{\linewidth}
    \epsfig{file=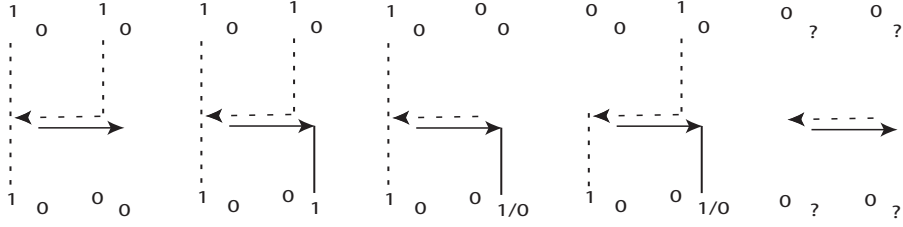, clip=,
          width=\duneu, height=3.0cm}
    \caption{\footnotesize The resampling mechanism 
       and the coalescent mechanism
       satisfy~\eqref{eq:g_bestimmt_f}}%
     \label{f:c_bestimmt_bd}
  \end{minipage}
  \end{beidschub}
  \end{figure}
\end{example}

The following proposition shows that two processes are dual in the sense
of~\eqref{eq:wedge_duality} if their defining basic mechanisms are dual
(cf.\ Definition~\ref{dualfunc}).
The proofs of both Proposition~\ref{dualcrit} and
Proposition~\ref{basicduality}
follow similar ideas as in~\cite{Gri79}.
\begin{proposition}\label{dualcrit}
  Let $m\in\N$ and let
  the processes $(X_t^N)_{t\geq0}$ and $(Y_t^N)_{t\geq0}$
  be defined by basic mechanisms
  $(f_1,...,f_m)$ and $(g_1,...$,$g_m)$, respectively.
  Suppose that the Poisson
  processes associated with $k\in \{1,\ldots,m\}$
  in $(X_t^N)_{t\geq0}$ and in $(Y_t^N)_{t\geq0}$
  have the same rate parameter $\ld_k\geq 0$.
  If $f_k$ and $g_k$ are dual for every $k=1,\ldots,m$,
  then $(X_t^N)_{t\geq0}$ and $(Y_t^N)_{t\geq0}$ satisfy the
  duality relation~\eqref{eq:wedge_duality}.
\end{proposition}
\begin{proof}
  Fix $T>0$ and initial values $X_0^N,Y_0^N\in\{0,1\}^N$.
  Assume for simplicity that $m=1$ and let $f:=f_1$, $g:=g_1$.
  Define the process $\rub{\hat{Y}_{t}^N}_{0\leq t\leq T}$ in
  backward time in the following way. Reverse all arrows in the graphical
  representation of $(X_t^N)_{t\geq0}$. At (forward) time $T$, start with a type
  configuration given by $\hat{Y}_0^N:=Y_0^N$.
  Now proceed until (forward) time $0$:
  Whenever you encounter an arrow, change the configuration according to $g$.
  Recall that the direction of the arrow indicates the order of the involved
  individuals.
  We show that the processes $(X_t^N)_{t\geq0}$ and
  $(\hat{Y}_t^N)_{0\leq t\leq T}$ satisfy
  \begin{equation}\label{eq:sanchez}
    X_0^N \wedge \hat{Y}_T^N = \uline{0}
    \,\Longleftrightarrow\,
    X_T^N \wedge \hat{Y}_0^N = \uline{0}
    \quad\fa X_0^N,\hat{Y}_0^N\in\{0,1\}^N,
  \end{equation}
  for every realisation.
  We prove the implication ``$\!\implies\!$'' by contradiction.
  Hence, assume that for some
  initial configuration there is a (random) time $t\in[0,T]$ such that
  \begin{equation}
    X_0^N \wedge \hat{Y}_T^N = \uline{0} \text{ and }
    X_t^N \wedge \hat{Y}_{T-t}^N \neq \uline{0}.
  \end{equation}
  There are only finitely many arrows until time $T$ and no two arrows occur
  at the same time almost surely.
  Hence, there is a first time $\tau$ such that the
  processes are disjoint before this time but not after this time.
  The arrow at time $\tau$ points from $i$ to $j$, say.
  Denote by $(x^-_i,x_j^{-})\in\{0,1\}^2$ and $(x_i^{+},x_j^{+})$ the types
  of the pair $(i,j)\in\{1,...,N\}^2$  
  according to the process $(X_t^N)_{t\geq0}$
  immediately before and after forward time $\tau$, respectively.
  By the definition of the process, we then have $f(x_i^{-},x_j^{-})=(x_i^{+},x_j^{+})$.
  Furthermore, denote by $(y_j^{-},y_i^{-})$ the types
  of the pair $(j,i)$ according to  $(Y_t^N)_{t\geq0}$
  immediately before backward time $T-\tau$.
  We have chosen $\tau,i,j$ such that
  \begin{equation}
    (x_i^{-},x_j^{-})\wedge \rub{g(y_j^{-},y_i^{-})}^\dagger=(0,0)
    \quad\text{ and }\quad
    (x_i^{+},x_j^{+})\wedge (y_i^{-},y_j^{-})\neq(0,0).
  \end{equation}
  However, this contradicts the duality of $f$ and $g$.
  The proof of the other implication is analogous.

  It remains to prove that $Y_T^N$ and $\hat{Y}_T^N$ are equal
  in distribution.
  The assertion then follows from
  \begin{equation}
    \P\eckb{X_0^N\wedge Y_T^N=\uline{0}}=
    \P\eckb{X_0^N\wedge \hat{Y}_T^N=\uline{0}}
    \stackrel{\eqref{eq:sanchez}}{=}
    \P\eckb{X_T^N\wedge \hat{Y}_0^N=\uline{0}}=
    \P\eckb{X_T^N\wedge Y_0^N=\uline{0}}.
  \end{equation}
  If a Poisson process is conditioned on its value at some fixed time $T > 0$, then the
  time points are uniformly distributed over the interval $[0,T]$.
  The uniform distribution is invariant under time reversal.
  In addition, the Poisson processes of $(Y_t^N)_{t\geq0}$ 
  nd $(X_t^N)_{t\geq0}$ have the same rate parameter.
  Thus, $(Y_t^N)_{0\leq t\leq T}$ and $(\hat{Y}_t^N)_{0\leq t\leq T}$
  have the same one-dimensional distributions.
\end{proof}

%
%
We will now give a list of those maps $f:\{0,1\}^2 \rightarrow \{0,1\}^2$
for which there exists a dual basic mechanism (see Definition \ref{dualfunc}).
The maps $f$ and $g$  in every row of the following table are dual to each other.
As in Example~\ref{ex:1}, it is elementary to check this.

\begin{center}
\begin{tabular}{l|c|c|c|c||c|c|c|c|}
$N^o$ & $f(0,0)$ & $f(0,1)$ & $f(1,0)$ & $f(1,1)$ & $g(0,0)$ & $g(0,1)$ & $g(1,0)$ & $g(1,1)$ \\ \hline \hline

i) & (0,0) & (0,0) & (1,1) & (1,1) & (0,0) & (0,1) & (0,1) & (0,1) \\ \hline
ii) & (0,0) & (0,1) & (1,1) & (1,1) & (0,0) & (0,1) & (1,1) & (1,1) \\ \hline
iii) & (0,0) & (0,0) & (0,1) & (0,1) & (0,0) & (0,0) & (0,1) & (0,1) \\ \hline
iv) & (0,0) & (0,1) & (1,0) & (1,1) & (0,0) & (0,1) & (1,0) & (1,1) \\ \hline
v) & (0,0) & (1,1) & (1,1) & (1,1) & (0,0) & (1,1) & (1,1) & (1,1) \\ \hline
vi) & (0,0) & (0,0) & (0,0) & (0,0) & (0,0) & (0,0) & (0,0) & (0,0) \\ \hline
\end{tabular}
\end{center}

\medskip\noindent
Check that the pair $(f,g)$ is dual if and only if the pair
$(f^\dagger,g^\dagger)$ is dual
where $f^\dagger(x):=(f(x^\dagger))^\dagger$.
Furthermore,
the pair $(f,g)$ is dual if and only if
$(\hat{f},\hat{g}^\dagger)$ is dual
where $\hat{f}(x):=f(x^\dagger)$
and $\hat{g}^\dagger(x)=(g(x))^\dagger$ for $x\in\{0,1\}^2$.
Thus, for each of the listed dual pairs $(f,g)$, 
the pairs  $(f^\dagger,g^\dagger)$,
$(\hat{f},\hat{g}^\dagger)$
and
$(\hat{f}^\dagger,\hat{g})$ are also dual.
Modulo this relation, the listing of dual basic mechanisms is complete.
The proof of this assertion is elementary but somewhat tedious and is 
deferred to the appendix.

Of particular interest are the dualities in i)-iii).
The first of these is the duality between the resampling mechanism
and the coalescent mechanism,
which we already encountered in Example~\ref{ex:1}.
The duality in ii)
is the self-duality of the \textbf{pure birth mechanism}
\begin{equation}\label{pbm}
  f^B:\{0,1\}^2 \rightarrow \{0,1\}^2, \, (1,0)\mapsto(1,1) \text{ and }
  x\mapsto x \; \forall\, x\in\{(0,0),(0,1),(1,1)\}
\end{equation}
and iii) is the self-duality of the \textbf{death/coalescent mechanism} 
\begin{equation}\label{dcm}
  f^{DC}:\{0,1\}^2 \rightarrow \{0,1\}^2, \, (1,\cdot) \mapsto (0,1) \text{ and }
    (0,\cdot) \mapsto (0,0).
\end{equation}
We are only interested in the effect of a basic mechanism on the total
number of individuals of type $1$. The identity map in iv) does not
change the number of individuals of type $1$ in the configuration.
The effect of v) and vi) on the number of individuals of type $1$ is similar to
the effect of ii) and iii), respectively.
Furthermore, both $f^\dagger$ and $\hat{f}$ have the same effect on the
number of individuals of type $1$ as $f$.

Closing this section,
we define processes which satisfy the duality
relation~\eqref{eq:wedge_duality}.
These processes will play a major role in deriving the
dualities \eqref{eq:braco_resem} and \eqref{eq:Fellog_dual} in Section 4.
For $u,e,\gamma,\beta \geq0$, 
let $(X_t^N)_{t\geq0}=(X_t^{N,(u,e, \gamma,\beta)})_{t\geq0}$
be the process on $\{0,1\}^N$ with the
following transition rates (of independent Poisson processes):
\begin{itemize}
  \item With rate $\tfrac{u}{N}$, the pure birth mechanism
     $f^B$ occurs (cf.\eqref{pbm}).
  \item With rate $\tfrac{e}{N}$, the death/coalescent mechanism
     $f^{DC}$ occurs (cf. \eqref{dcm}).
  \item With rate $\tfrac{\gamma}{N}$, the coalescent mechanism
     $f^C$ occurs (cf.\ \eqref{eq:fC}).
  \item With rate $\tfrac{\beta}{N}$, the resampling mechanism
     $f^{R}$ occurs (cf.\ \eqref{eq:fC}).
\end{itemize}
Together with an initial configuration, this defines the process.
The process $(X_t^{N,(u,e, \gamma,\beta)})_{t\geq0}$
is defined by the basic
mechanisms $(f^B,f^{DC},f^C,f^R)$,
and  the process $(X_t^{N,(u,e,\beta ,\gamma)})_{t\geq0}$ is
defined by the basic mechanisms  $(f^B,f^{DC},f^R,f^C)$.
Proposition~\ref{dualcrit} then yields the following corollary.
\begin{corollary}    \label{c:satisfy_wedge}
  Let $u,e,\gamma,\beta\geq0$.
  The two processes  $(X_t^{N,(u,e, \gamma,\beta)})_{t\geq0}$
  and  $(X_t^{N,(u,e,\beta ,\gamma)})_{t\geq0}$
  satisfy the duality relation~\eqref{eq:wedge_duality}.
\end  {corollary}
%
%
%
\section{Prototype duality}%
\label{sec:prototype_duality}
In this section, we derive the prototype duality~\eqref{dadudi}
from~\eqref{eq:wedge_duality}. The main idea for this is
to integrate equation~\eqref{eq:wedge_duality} in the variables $x^N$ and
$y^N$ with respect to a suitable measure.
Furthermore, we will exploit the fact that drawing from an urn with
replacement and without replacement, respectively, is almost surely
the same if the urn contains infinitely many balls.
%
%
\enlargethispage{0.2cm}
\begin{proposition}\label{basicduality}
  Let $(X_t^N)_{t\geq0}$ and $(Y_t^N)_{t\geq0}$ be processes
  with state space $\{0,1\}^N$, $N\geq 1$.
  Assume that $(X_t^N)_{t\geq0}$ and $(Y_t^N)_{t\geq0}$ satisfy the duality
  relation~\eqref{eq:wedge_duality}.
  Choose $n,k\in\{0,...,N\}$
  which may depend on $N$.
  Define $\mu_n^N(x^N):=\binom{N}{n}^{-1}\1_{\abs{x^N}=n}$ for
  every $x^N\in\{0,1\}^N$ where $\abs{x^N}=\sum_{i=1}^N x_i^N$ is the
  total number of individuals of type $1$.
  Assume $\Law{X_0^N}=\mu_n^N$ and $\Law{Y_0^N}=\mu_k^N$.
  Suppose that the process $(X_t^N)_{t\geq0}$ satisfies
  \begin{equation}  \label{a1}
    \frac{n}{N}\to0\quad\text{and}\quad
    \frac{\E\eckb{\absb{X_{t_N}^N}}}{N}\,{\longrightarrow}\, 0\qqasN,
  \end  {equation}
  where $t_N\geq0$.
  Then
  \begin{equation}  \label{dualrel}
    \limN \E\left[\ruB{1-\tfrac{k}{N}}^{\absb{X_{t_N}^N}}\right] 
    =\limN \E\left[\ruB{1-\tfrac{\absb{Y_{t_N}^N}}{N}}^n\right] 
  \end{equation}
  under the assumption that the limits exist.
\end{proposition}

\begin{proof}
A central idea of the proof is to make use of the well known fact that the
hypergeometric distribution $\text{Hyp}(N,R,l)$, $R,l\in\{0,...,N\}$,
can be
approximated by the binomial distribution $\text{B}(l,\tfrac{R}{N})$
as $N\to\infty$
provided that $l$ is sufficiently small compared to $N$.
In fact, by Theorem 4 of~\cite{DF80},
\begin{equation}
  \absB{\text{B}(l,\tfrac{R}{N})\eckb{\{0\}}-\text{Hyp}(N,R,l)\eckb{\{0\}}}
  \leq d_{TV}\rubb{\text{B}\rub{l,\tfrac{R}{N}},\text{Hyp}(N,R,l)}
  \leq \frac{4 \cdot l}{N} \quad\fa R,l\leq N,
\end{equation}
where $d_{TV}$ is the total variation distance.
By assumption~\eqref{a1}, we have (with $R:=k, l:=\absb{X_{t_N}^N}$)
\begin{equation}\label{a4}
  \E\eckbb{\ruB{1-\tfrac{k}{N}}^{\absb{X_{t_N}^N}}}
  = \E\eckbb{\text{B}\ruB{\absb{X_{t_N}^N},\tfrac{k}{N}}\eckb{\{0\}}}
  = \E\eckB{\text{Hyp}\rub{N,k,\absb{X_{t_N}^N}}\eckb{\{0\}}} + o(1)
\end{equation}
as $N\to\infty$. Similarly, we have (with $R:=\absb{Y_{t_N}^N}$, $l:=n$)
\begin{equation}\label{a3}
  \E\eckbb{\ruB{1-\tfrac{\absb{Y_{t_N}^N}}{N}}^n}
  = \E\eckbb{\text{B}\ruB{n,\tfrac{\absb{Y_{t_N}^N}}{N}}\eckb{\{0\}}}
  = \E\eckB{\text{Hyp}\rub{N,\absb{Y_{t_N}^N},n}\eckb{\{0\}}} +o(1)
\end{equation}
as $N\to\infty$.
%
By definition of the hypergeometric distribution, we get
\begin{equation}\label{a5}
  \text{Hyp}\rub{N,\absb{Y_t^N},n}\eckb{\{0\}}
  = \binom{N}{n}^{-1} \sum_{x^N: \, |x^N|= n}
   \mathbbm{1}_{\{ x^N\wedge \, Y_t^N = \uline{0}\}}
  =  \mu_n^N\eckb{x^N\colon x^N \wedge Y_t^N = \uline{0}}.
\end{equation}
By the same argument, we also obtain
\begin{equation}\label{a6}
  \text{Hyp}\rub{N,k,\absb{X_t^N}}\eckb{\{0\}}
  =\text{Hyp}\rub{N,\absb{X_t^N},k}\eckb{\{0\}}
  = \mu_k^N\eckb{y^N\colon X_t^N \wedge y^N = \uline{0}}.
\end{equation}
We denote by $\P^{x^N}$ the law of the process $(X_t^N)_{t\geq0}$ started
in the fixed initial configuration $x^N\in\{0,1\}^N$.
Starting from the left-hand side of~\eqref{dualrel},
the above considerations yield
\begin{equation}\begin{split}
  \E& \eckB{ \rub{1-\frac{k}{N}}^{\absb{X_{t_N}^N}} } +o(1)
  \stackrel{\eqref{a4}}{=}
     \E\eckB{ \text{Hyp}(N,k,\absb{X_{t_N}^N})\eckb{\{0\}}}\\
  \etStackrel{\eqref{a6}}{=}
    \int\E^{x^N}\eckB{\mu_k^N\eckb{X_{t_N}^N \wedge y^N=\ul{0}}}
       \mu_n^N(\dup x^N)\\
  \etStackrel{\eqref{eq:wedge_duality}}{=}
    \int \int
     \P^{y^N}\eckB{x^N \wedge Y_{t_N}^N=\ul{0}}
    \mu_k^N(\dup y^N)\,\mu_n^N(\dup x^N)
  \stackrel{}{=}
    \E\eckB{ \mu_n^N\eckb{ x^N \wedge Y_{t_N}^N = \ul{0}}}\\
  \etStackrel{\eqref{a5}}{=}
    \E\eckB{\text{Hyp}(N,\absb{Y_{t_N}^N},n)\eckb{\{0\}}} 
  \stackrel{\eqref{a3}}{=} 
    \E\eckbb{\ruB{1-\frac{\absb{Y_{t_N}^N}}{N}}^n } + o(1),
\end{split}
\end{equation}
which proves the assertion.
\end{proof}

\section{Various scalings}%
\label{sec:different_scalings}

Recall the definition of the process
$(X_t^{N,(u,e,\gamma,\beta)})_{t\geq0}$
from the end of Section~\ref{sec:dual_basic_mechanisms}.
Define $X_t^N:=X_t^{N,(u,e,\gamma,\beta)}$ and
$Y_t^N:=X_t^{N,(u,e,\beta,\gamma)}$ for $t\geq0$ and $N\in\N$.
Notice that the Poisson process attached to the resampling mechanism
in the process $(Y_t^N)_{t\geq0}$ has rate $\gamma$.
By Corollary~\ref{c:satisfy_wedge}, the two processes
$(X_t^N)_{t\geq0}$ and
$(Y_t^N)_{t\geq0}$ satisfy the duality
relation~\eqref{eq:wedge_duality}.
Let $\Law{X_0^N}=\mu_n^N$ and $\Law{Y_0^N}=\mu_k^N$ for some
$n,k\in\N$ to be chosen later, where $\mu_n^N$ is defined in
Proposition~\ref{basicduality}.
In order to apply Proposition~\ref{basicduality},
we essentially have to prove existence of the limits in~\eqref{dualrel}.
Depending on the scaling, this will result in the
moment duality~\eqref{eq:braco_resem}
of a resampling-selection model with
a branching-coalescing particle process
and
in the Laplace duality~\eqref{eq:Fellog_dual} of the logistic Feller diffusion with
another logistic Feller diffusion,
respectively.
Both dualities
could be derived
simultaneously.
However, in order to keep things
simple, we consider the two cases separately.
%
%
%
%
%
\begin{theorem}\label{GWdual}
Assume that $b,c,d\geq0$.
Denote by $(X_t)_{t\geq0}$ and $(Y_t)_{t\geq0}$
the $(1,b,c,d)$-braco-process
and the $(1,b,c,d)$-resem-pro\-cess, respectively.
The initial values
are $X_0=n\in\N_0$ and $Y_0=y\in[0,1]$.
Then 
\begin{equation}
  \E^n\eckb{(1-y)^{X_t}} =\E^y\eckb{(1-Y_t)^n},\qquad t\geq0.
\end{equation}
\end{theorem}
\begin{remark} \label{r:GWdual}
  In the special case $b=0=d$ and $c>0$, this is the moment
  duality of the
  Fisher-Wright diffusion with Kingman's coalescent.
  Furthermore, choosing $c=0$ and $b,d>0$ results in the
  moment duality of the
  Galton-Watson process with a deterministic process.
\end{remark}
\begin{proof}
  Choose $u,e,\beta\geq0$ and $\gamma=\gamma(N)$ such that
  $b=u+\beta$, $d=e+\beta$ and $\gamma/N\to c$ as $N\to\infty$.
  In the first step, we prove that the process $(\abs{X_t^N})_{t\geq0}$
  of the total number of individuals of type $1$ 
  converges weakly to $(X_t)_{t\geq0}$.
  The total number of individuals of type $1$
  increases by one if a ``birth event'' occurs ($f^B$ or $f^{R}$)
  and if the type configuration of the respective ordered pair
  of individuals is $(1,0)$.
  If the total number of individuals of type $1$ is equal
  to $k$, then the probability
  of the type configuration of a randomly chosen ordered pair to
  be $(1,0)$ is $\tfrac{k}{N} \tfrac{N-k}{N-1}$.
  The number of Poisson processes associated with a fixed basic
  mechanism is $N(N-1)$. Thus, the process of the total
  number of individuals of type $1$
  has the following transition rates:
  \begin{equation} \label{eq:transition_rates}
  \begin{array}{ll}
    k\rightarrow k+1: & \frac{u+\beta}{N} \cdot N(N-1)
      \cdot \frac{k}{N} \frac{N-k}{N-1}, \\
    k\rightarrow k-1: & \frac{e+\beta}{N} \cdot N(N-1)
      \cdot \frac{N-k}{N} \frac{k}{N-1} + \frac{e+\gamma}{N}
      \cdot N(N-1) \cdot \frac{k}{N} \frac{k-1}{N-1},
  \end{array}   
  \end{equation}
  where  $k\in \N_0$.
  Notice that the coalescent mechanism produces the quadratic term
  $k(k-1)$ because the probability
  of the type configuration of a randomly chosen ordered pair to
  be $(1,1)$ is $\tfrac{k}{N} \tfrac{k-1}{N-1}$
  if there are $k$ individuals of type $1$.
  The transition rates determine
  the generator $\Gen^N=\Gen^{N,(u,e,\gamma,\beta)}$
  of $(\abs{X_t^N})_{t\geq0}$, namely
  \begin{equation}\begin{split}   \label{CGN}
    \Gen^{N} f(k) = & \frac{u+\beta}{N} \cdot k(N-k)
      \cdot \big(f(k+1)-f(k)\big)\\ 
    &+ \frac{e+\beta}{N} \cdot k(N-k) \cdot \big(f(k-1) -f(k)\big)\\
    &+ \frac{e+\gamma}{N} \cdot k(k-1) \cdot \big(f(k-1)-f(k)\big),
    \quad k\in\{0,\ldots,N\},
  \end{split} \end{equation} 
  for $f\colon\{0,\dots,N\}\to\R$.
  The $(1,u+\beta,c,e+\beta)$-braco-process
  $(X_t)_{t\geq0}$ is the unique solution of the martingale problem 
  for $\Gen$ (see~\cite{AS05}) where
  \begin{equation}
    \Gen f(k) := (u+\beta) k \, \rub{f(k+1)-f(k)}
          +\rub{(e+\beta) +c(k-1)}k\,\rub{f(k-1) - f(k)}, \quad k\in \N_0,
  \end{equation}
  for $f\colon\N_0\to\R$ with finite support.
  Letting $N\to \infty$, we see that
  \begin{equation}
    \Gen^{N} f(k) \lra\Gen f(k)\qqasN,\quad k\in\N_0,
  \end{equation}
  for $f\colon\N_0\to\R$ with finite support.
  We aim at using Lemma~\ref{weakconv} which is given below
  (with $E_N=\{0,\ldots,N\}$ and $E=\N_0$), to infer from (\ref{CGN})
  the weak convergence of the corresponding Markov processes.
  A coupling argument shows that $(\abs{X_t^{N}})_{t\geq0}$
  is dominated by
  $(Z_t^N)_{t\geq0}:= (\abs{X_t^{N,(u,0,0,\beta)}})_{t\geq0}$.
  The process $(Z_t^N)_{t\geq0}$ solves the martingale problem
  for $\Gen^{N,(u,0,0,\beta)}$.
  Thus, we obtain
  \begin{equation}
    Z^N_t-Z^N_0=\int_0^t\Gen^{N,(u,0,0,\beta)} Z^N_s\,ds +C^N_t
     =\int_0^t u Z^N_s\tfrac{N-Z_s^N}{N}\,ds +C^N_t
  \end{equation}
  where $(C_t^N)_{t\geq0}$ is a martingale.
  Hence, $(Z_t^N)_{t\geq0}$ is a submartingale.
  Taking expectations, Gronwall's inequality implies
  \begin{equation}  \label{eq:Z_first_moment}
    \E [Z_t^N]\leq \E [Z_0^N] e^{ut},  \qquad \fa t\geq0.
  \end{equation}
  Let $S_N=T_N=1$, $s_N=u$ and recall $\abs{X_0^N}=n$.
  With this, the assumptions of Lemma~\ref{weakconv} are satisfied.
  Thus, Lemma~\ref{weakconv} implies that
  $(\abs{X_t^N})_{t\geq0}$
  converges weakly to $(X_t)_{t\geq0}$ as $N\to\infty$.
  Let $k=k_N\in\{0,...,N\}$ be such that $k/N\to y$ as $N\to\infty$.
  For every $\bar{n}\in\N$,
  $(1-\tfrac{k}{N})^n$ converges uniformly in $n\leq\bar{n}$
  to $(1-y)^n$ as $N\to\infty$.
  In general, if the sequence $(\tilde{X}_n)_{n\in\N}$ of random
  variables with complete and separable state space converges weakly
  to $\tilde{X}$ and if the sequence $(f_n)_{n\in\N}$,
  $f_n\in{C}_b$,
  converges uniformly on compact sets to $f\in{C}_b$,
  then $\E [f_n(\tilde{X}_n)]\to\E [f(\tilde{X})]$ as $n\to\infty$.
  Hence,
  \begin{equation}  \label{eq:erster_limes}
    \E^n\eckB{(1-y)^{X_t}}=\lim_{N \to \infty}
     \E \eckbb{\Big(1-\frac{k}{N}\Big)^{\absb{X^{N}_t}}}.
  \end{equation}

  The next step is to prove that the rescaled processes
  $\ru{\abs{Y_t^{N}}/{N}}_{t\geq0}$
  converge weakly to $(Y_t)_{t\geq0}$ as $N\to\infty$.
  The generator of $\ru{\abs{Y_t^{N}}/{N}}_{t\geq0}$ is given by
  \begin{equation}\begin{split}   \label{eq:genesis}
    \Gen^{N,(u,e,\beta,\gamma)} f\rub{\tfrac{k}{N}}=
     &\gamma k \,\frac{N-k}{N}\Big(f\big(\tfrac{k+1}{N}\big)
         +f\big(\tfrac{k-1}{N}\big) - 2f\big(\tfrac{k}{N}\big)\Big)\\  
     &+u k \,\frac{N-k}{N}\,\rubb{f\rub{\tfrac{k+1}{N}} -
         f\rub{\tfrac{k}{N}} }
      + e k \,\frac{N-k}{N}\Big(f\big(\tfrac{k-1}{N}\big) -
         f\big(\tfrac{k}{N}\big)\Big)\\  
     &+ \tfrac{e+\beta}{N} \, k(k-1) \, \Big(  f\big(\tfrac{k-1}{N}\big) -
         f\big(\tfrac{k}{N}\big) \Big),
      \qquad k\in\{0,...,N\},
  \end{split} \end{equation}
  for $f \in C_c^2([0,1])$.
  Choose $k = k_N\leq N$ such that
  $\tfrac{k}{N} \rightarrow y \in [0,1]$ as $N\to\infty$.
  Notice that
  \begin{equation}
    N^2\mal\ruB{f\rub{\tfrac{k+1}{N}}+f\rub{\tfrac{k-1}{N}}
      -2f\rub{\tfrac{k}{N}}} \to f^{''}(y)
    \qasN.
  \end{equation}
  As $N\to\infty$,
  the right-hand side of~\eqref{eq:genesis} converges to
  \begin{equation}\begin{split}
    \lefteqn{c y(1-y)\cdot f^{''}(y)
      +(u-e) y (1-y) \cdot f^{'}(y)
      - (e+\beta) y^2\cdot  f^{'}(y)}\\
    &= (u-e) y \cdot f^{'}(y) - (u+\beta) y^2 \cdot f^{'}(y)
      +c y(1-y)\cdot f^{''}(y)=:\Gen f(y)
  \end{split}
  \end{equation}
  for every $f\in C_c^2([0,1])$.
  Athreya and Swart~\cite{AS05} show that
  the $(1,b,c,d)$-resem-process $(Y_t)_{t\geq0}$ solves the martingale
  problem for $\Gen$ and that this solution is unique.
  Let
  $E_N=\{0,1,\ldots,N\}$,
  $E=[0,1]$, $Z_t^N:=\abs{X_t^{N,(u,0,0,\gamma)}}$, $S_N=N$ and $T_N=1$.
  With this, the assumptions of Lemma~\ref{weakconv}
  are satisfied and we conclude that
  $(\abs{Y_t^{N}}/N)_{t\geq0}$ converges weakly to $(Y_t)_{t\geq0}$.
  It follows that, for $k=k_N \in \{0,...,N\}$
  with $k/N\to y$,
  \begin{equation}  \label{eq:zweiter_limes}
    \limN \E\left[\ruB{1-\frac{\absb{Y_{t}^N}}{N}}^n\right] 
    =\E^y\eckB{(1-Y_t)^{n}}.
  \end  {equation}
  This proves existence of the limits in~\eqref{dualrel}
  with $t_N:=t$.
  Inequality~\eqref{eq:Z_first_moment} and
  $\abs{X_0^{N}}=n<<N$ imply condition~\eqref{a1}.
  Thus, Proposition~\ref{basicduality}
  establishes equation~\eqref{dualrel}.
  The assertion follows from equations~\eqref{eq:erster_limes},
  \eqref{dualrel} and~\eqref{eq:zweiter_limes}.
\end{proof}


%
%
Next, we derive the Laplace duality of a logistic Feller diffusion
with another logistic Feller diffusion.
Recall that the logistic Feller diffusion with parameters
$(\al,\gamma,\beta)$ solves equation~\eqref{eq:logistic_Feller}.
\begin{theorem}\label{Felldual}
  Suppose that $\al,\gamma,\beta\geq0$, $r>0$ and $X_0=x\geq0,Y_0=y\geq0$.
  Let $(X_t)_{t\geq0}$ and $(Y_t)_{t\geq0}$ be logistic Feller
  diffusions with parameters
  $(\al,\gamma,\beta)$ and $(\al,r \beta,\gamma/r)$, respectively.
  Then
  \begin{equation}
    \E^x \eckb{e^{-r X_t\cdot y}} =\E^y \eckb{ e^{-r x\cdot Y_t}}
  \end{equation}
  for all $t\geq0$.
\end{theorem}
\begin{remarkb} \label{r:Felldual}
  \begin{enumerate}[(a)]
    \item  For $\beta,\gamma>0$ and $r=\gamma/\beta$, 
      Theorem~\ref{Felldual} yields the self-duality of the logistic Feller diffusion.
    \item  For  $\al=0$, $\gamma=0$, $r=1$ and $\beta > 0$,  Theorem~\ref{Felldual}
      specialises to the Laplace duality of Feller's branching diffusion.
  \end  {enumerate}
\end{remarkb}
\begin{proof}
  Choose $u=u_N\geq0$ and $e=e_N\geq0$ such that $(u-e)\sqrt{N}\to\al$
  as $N\to\infty$.
  We prove that the rescaled process
  $(\abs{Y_{t\sqrt{N}}^N}/(r\sqrt{N}))_{t\geq0}$
  converges weakly to $(Y_t)_{t\geq0}$ as $N\to\infty$.
  The generator of the rescaled process 
  is given by  (cf.\ \eqref{eq:genesis})
  \begin{equation} \begin{split}   \label{eq:adidas}
    \sqrt{N}\Gen^{N} f\rub{\tfrac{k}{r\sqrt{N}}}
    &= \sqrt{N} \cdot \gamma \cdot k\frac{(N-k)}{N} \cdot
       \Big( f\big(\tfrac{k+1}{r\sqrt{N}}\big) +
       f\big( \tfrac{k-1}{r\sqrt{N}} \big)-2f\rub{\tfrac{k}{r\sqrt{N}}}\Big)\\
    &\quad + \sqrt{N}u_N \cdot k\frac{(N-k)}{N} \cdot
      \ruB{ f\big( \tfrac{k+1}{r\sqrt{N}} \big) -f\rub{\tfrac{k}{r\sqrt{N}}}}\\
    &\quad +  \sqrt{N}e_N \cdot k\frac{(N-k)}{N} \cdot
      \Big( f\big( \tfrac{k-1}{r\sqrt{N}} \big) -
      f\big( \tfrac{k}{r\sqrt{N}} \big) \Big)\\
    &\quad + \sqrt{N} \cdot (e_N + \beta)\cdot 
      \frac{k(k-1)}{r^2N}\,r^2  \cdot
      \frac{r\sqrt{N}}{r\sqrt{N}}\,\Big( f\big( \tfrac{k-1}{r\sqrt{N}} \big) -
      f\big( \tfrac{k}{r\sqrt{N}} \big) \Big),
  \end{split} \end{equation}
  for $k\in\{0,\ldots,N\}$ and for $f \in C^2_c([0,\infty))$.
  Let $k=k(N)\in\{0,\ldots,N\}$ be such that
  $k/(r\sqrt{N})\to y$. Letting $N\to\infty$, the right-hand side
  converges to
  \begin{equation}  \label{eq:puma}
    \tfrac{\gamma }{r}\,y\cdot f^{''}(y)+\al\,y\cdot f^{'}(y)
       -\beta r\, y^2\cdot f^{'}(y)=:\Gen f(y)
  \end{equation}
  for every $f \in C^2_c([0,\infty))$.
  Notice that the quadratic term $y^2$ originates in  the
  quadratic term $k(k-1)$.
  Hutzenthaler and Wakolbinger~\cite{HW07} prove that $(Y_t)_{t\geq0}$
  is the unique solution of the martingale problem for $\Gen$.
  Let $\abs{Y_0^N}=k=k(N)$ be such that 
  $k/(r\sqrt{N})\to y\in[0,1]$ as $N\to\infty$ and define $Z_0^N:=k$.
  As before, $(Z_t^N)_{t\geq0}:=(\abs{X_t^{N,(u,0,0,\gamma)}})_{t\geq0}$
  is a submartingale which dominates $(Y_t^N)_{t\geq0}$ and which satisfies
  \begin{equation} \label{eq:Z_Laplace}
    \sup_N\frac{1}{r\sqrt{N}}\E [Z_{t\sqrt{N}}^N]
    \leq \sup_N\frac{1}{r\sqrt{N}}\E [Z_0^N] e^{u_Nt\sqrt{N}}<\infty, 
      \qquad \fa t\geq0.
  \end{equation}
  Let $E_N:=\{0,\ldots,N\}$,
  $E:=[0,\infty)$, $s_N:=u_N$, $S_N:=r\sqrt{N}$ and $T_N:=\sqrt{N}$.
  The assumptions of Lemma~\ref{weakconv}
  are satisfied and we conclude that
  $(\abs{Y_{t\sqrt{N}}^N}/(r\sqrt{N}))_{t\geq0}$
  converges weakly to $(Y_t)_{t\geq0}$.
  This also proves that 
  $(\abs{X_{t\sqrt{N}}^N}/\sqrt{N})_{t\geq0}$
  converges weakly to $(X_t)_{t\geq0}$ if
  $\abs{X_0^N}=n=n(N)$ is such that $n/\sqrt{N}\to x$ as $N\to\infty$. It is not hard to see that, for every $\tilde{z}\geq0$,
  \begin{equation}
    \rub{1-r\tfrac{k/(r\sqrt{N})}{\sqrt{N}}}^{\sqrt{N}z}\lra e^{-rzy}
    \quad\text{ and }\quad
    \rub{1-r\tfrac{z}{\sqrt{N}}}^{\sqrt{N}\tfrac{n}{\sqrt{N}}}\lra e^{-rxz}
    \qqasN
  \end{equation}
  uniformly in $0\leq z\leq\tilde{z}$. Together with the weak convergence
  of the rescaled processes, this implies
  \begin{equation} \label{eq:erster_Laplace}
    \E^x \eckb{e^{-r X_t\cdot y}} =
    \limN \E^n \left[ \Big(1-r\tfrac{k/(r\sqrt{N})}{\sqrt{N}}\Big)^{\sqrt{N}
      \cdot X^N_{t\sqrt{N}}/\sqrt{N}} \right] 
  \end{equation}
  and
  \begin{equation} \label{eq:zweiter_Laplace}
    \limN\E^k\left[\Big(1-r\tfrac{Y^N_{t\sqrt{N}}
         /(r\sqrt{N})}{\sqrt{N}}\Big)^n \right]
    =\E^y \eckb{ e^{-r x\cdot Y_t}}
  \end{equation}
  for $t\geq0$.
  This proves existence of the limits in~\eqref{dualrel}
  with $t_N:=t\sqrt{N}$.
  Inequality~\eqref{eq:Z_Laplace} and
  $\abs{X_0^{N}}=n<<N$ imply condition~\eqref{a1}.
  Thus, Proposition~\ref{basicduality}
  establishes equation~\eqref{dualrel}.
  The assertion follows from equations~\eqref{eq:erster_Laplace},
  \eqref{dualrel} and~\eqref{eq:zweiter_Laplace}.
\end{proof}
\begin{remark}   \label{r:pure_death_process}
  Assume $u=e=\gamma=\alpha=0$ and $r=1$ in the proof of
  Theorem~\ref{Felldual}. Then $\ru{\abs{Y_t^N}}_{t\geq0}$ is
  a pure death process on $\{1,...,N\}$ which jumps from $k$
  to $k-1$ at exponential rate $\tfrac{\beta}{N}k(k-1)$, $2\leq k\leq N$.
  Furthermore, $(Y_t)_{t\geq0}$ is a solution of~\eqref{eq:def:y_t}.
  We have just shown that the rescaled pure death process
  $\ru{\abs{Y_{t\sqrt{N}}^N}/\sqrt{N}}_{t\geq0}$
  converges weakly to $(Y_t)_{t\geq0}$ as $N\to\infty$.
\end  {remark}

\section{Weak convergence of processes}%
\label{sec:weak_convergence}
 
In the proofs of Theorem~\ref{GWdual} and Theorem~\ref{Felldual}, we have
established convergence of generators plus a domination principle.
In this section, we prove that this implies weak convergence of the
corresponding processes.
For the weak convergence of processes with  c\`adl\`ag paths,
let the topology on the set of  c\`adl\`ag paths be given
by the Skorohod topology (see \cite{EthKu}, Section 3.5).

\begin{lemma}\label{weakconv}
Let $E\subset\R_{\geq0}$ be closed.
Assume that
the martingale problem for $(\Gen,\nu)$ has at most one solution where
$\Gen\colon C_c^2(E)\to C_b(E)$ is a linear operator and
$\nu$ is a probability measure on $E$.
Furthermore, for $N\in\N$, let $E_N\subset \R_{\geq0}$ and let
$(Y_t^N)_{t\geq0}$ be an $E_N$-valued Markov process with c\`adl\`ag paths
and generator $\Gen^N$.
Let
$(S_N)_{N\in\N}$ and $(T_N)_{N\in\N}$
be sequences in\/ $\R_{>0}$ with $y^N/S_N\in E$ for all $y^N\in E_N$
and $N\in\N$.
Suppose that
\begin{equation}  \label{eq:generator_konv}
   y^N\in E_N, \limN \tfrac{y^N}{S_N}=y\in E\ \text{ implies }\  
     T_N\Gen^N f\rub{\tfrac{y^N}{S_N}}
     \rightarrow \Gen  f(y) \qqasN,
\end{equation}
for every $f \in C_c^2(E)$.
Assume that, for $N\in\N$, $(Y_t^N)_{t\geq0}$ is dominated by a
process $(Z_t^N)_{t\geq0}$, i.e.,
$Y_t^N \leq Z_t^N$ for all $t\geq0$ almost surely, which is
a submartingale satisfying $\E[ Z_t^N]\leq \E [Z_0^N] e^{ts_N}$
for all $t\geq0$ and some constant $s_N$.
In addition, suppose that \/$\tlimsupN s_N T_N < \infty$ and
\/$\tlimsupN \tfrac{\E[Z_0^N]}{S_N} < \infty$.
If $Y_0^N/S_N$ converges in distribution to $\nu$ as $N\to\infty$, then
\begin{equation}
  \LawB{\rub{{Y_{tT_N}^N}\big/{S_N}}_{t\geq0}}
  \wlim \LawB[\nu]{\rub{Y_t}_{t\geq0}}\qqasN
\end{equation}
where $(Y_t)_{t\geq0}$ is a solution of
the martingale problem $(\Gen,\nu)$
with initial distribution $\nu$.
\end{lemma}

\begin{proof}
We aim at applying Corollary 4.8.16 of Ethier and Kurtz \cite{EthKu}.
For this, define
\begin{equation}
  \tilde{E}_N:=\{\tfrac{y^N}{S_N}\colon y^N\in E_N\},
  \quad \tilde{\Gen}^N f(\tilde{y}^N)
    :=T_N\Gen^N f\rub{\tfrac{y^N}{S_N}}\Big|_{y^N=\tilde{y}^N S_N},
    \tilde{y}^N\in\tilde{E}^N,
\end{equation}
for $f\in C_c^2(E)$
and let $\eta_N\colon \tilde{E}_N\to E$ be the embedding function.
The process $\rub{Y_{tT_N}^N/S_N}_{t\geq0}$ has state space $\tilde{E}_N$
and generator $\tilde{\Gen}^N$.
Now we prove the compact containment condition, i.e.,
for fixed $\eps,t > 0$ we show
\begin{equation}
  \rub{\exists K > 0} \, \rub{\forall N \in \N} \;
  \P\eckB{ \sup_{s\leq t}
  \frac{Y_{sT_N}^N}{S_N} \leq K} \geq 1-\eps.
\end{equation}
Using $Y_t^N\leq Z_t^N$, $t\geq0$, and Doob's Submartingale Inequality,
we conclude for all $N\in\N$
\begin{equation}  \begin{split}
  \P\eckB{ \sup_{s\leq t} Y_{sT_N}^N \geq K S_N} & \leq
    \P\eckB{ \sup_{s\leq t} Z_{sT_N}^N \geq K S_N}
    \leq \frac{1}{K S_N} \E\eckb{Z^N_{tT_N}} \\
   & \leq \frac{1}{K} \sup_{N\in\N}\frac{\E\eckb{Z^N_0}}{S_N}
        \cdot \exp\rub{t \cdot \sup\limits_{N\in\N} (s_N T_N)}
    =:\frac{C}{K}.
\end{split}     \end{equation}
Thus, choosing $K := \tfrac{C}{\eps}$ completes the proof
of the compact containment condition.

It remains to verify condition (f) of Corollary 4.8.7 of~\cite{EthKu}.
Condition~\eqref{eq:generator_konv} implies that for
every $f\in C_c^2$ and every compact set $K\subset E$
\begin{equation} \label{Gen_unif}
  \sup_{y\in K\cap \tilde{E}_N}| \tilde{\Gen}^N f(y) - \Gen  f(y)|
     \rightarrow 0 \qqasN.
\end{equation}
Choose a sequence $K_N$ such that~\eqref{Gen_unif} still holds
with $K$ replaced by $K_N$. 
This together with the compact containment condition 
implies condition (f) of Corollary 4.8.7 of \cite{EthKu}
with $G_N:=K_N\cap \tilde{E}_N$ and $f_N:=f|_{\tilde{E}_N}$. 
Furthermore, notice that $C_c^2(E)$ is an algebra that separates points
and $E$ is complete and separable.
Now Corollary 4.8.16 of Ethier and Kurtz \cite{EthKu} implies
the assertion.
\end{proof}

\medskip
\noindent
\textbf{\large{Open Question:}} Athreya and Swart~\cite{AS05} %
prove a self-duality of the resem-process
given by~\eqref{eq:resem}.
We were not able to establish
a graphical representation for this duality.
Thus, the question whether our technique also works in this case
yet waits to be answered. 

\bigskip
\noindent
\textbf{\large{Acknowledgements:}}
We thank Achim Klenke and Anton Wakolbinger
for valuable discussions and many detailed remarks.
Also, we thank the referee for a number of very helpful suggestions.

\section*{Appendix}
\renewcommand{\thelemma}{A.\arabic{lemma}}
\setcounter{lemma}{0}

The aim here is to provide a complete list of dual basic ``mechanisms'', i.e. all combinations of maps $f,g: \{0,1\}^2 \rightarrow \{0,1\}^2$ fulfilling the conditions
\begin{eqnarray}
  (\forall \, x,y \in \{0,1\}^2) \;\;
       y \wedge f(x)^\dagger = (0,0) \,
  &\Longrightarrow&g(y) \wedge x^\dagger = (0,0)\quad
     \text{and}\label{eq:a:f_bestimmt_g}\\
  (\forall \, x,y \in \{0,1\}^2)\; \;
       x \wedge g(y)\dagger = (0,0) \,
  &\Longrightarrow& f(x) \wedge y^\dagger = (0,0),\label{eq:a:g_bestimmt_f}
\end{eqnarray}
where for $x=(x_1,x_2)\in \{0,1\}^2$ we define $x^\dagger := (x_2,x_1)$.\\
Define for a basic mechanism $f:\{0,1\}^2 \rightarrow \{0,1\}^2$ the maps $f^{\dagger}$ and $\hat{f}$ via $f^{\dagger}(x) := f(x^\dagger)^\dagger$ and $\hat{f}(x):= f(x^\dagger)$, leading to a third map $\hat{f}^{\dagger}$  with $\hat{f}^{\dagger}(x):= f^{\dagger}(x^\dagger) = \hat{f}(x^\dagger)^\dagger =  f(x)^\dagger$.  Then we get the following characterisation of duality between basic mechanisms $f,g$.
\begin{lemma}
Let $f,g$ be basic mechanisms, then
\begin{equation}\label{dbm}
f  \text{ and } g \text{ are dual}
 \Longleftrightarrow \ruB{(\forall x,y \in \{0,1\}^2)\; y \wedge f(x)= (0,0)  \, \Leftrightarrow \,  g^{\dagger}(y) \wedge x = (0,0)}.
\end{equation}
\end{lemma}
\begin{proof}
"$\Rightarrow$'': Let $f,g$ be dual, then (\ref{eq:a:f_bestimmt_g}) and (\ref{eq:a:g_bestimmt_f}) hold and thus we have for all $x,y \in \{0,1\}^2$:
\begin{equation*}
  y \wedge f(x) = (0,0)
\,\Longleftrightarrow \, y^\dagger \wedge f(x)^\dagger = (0,0)
\, \stackrel{(\ref{eq:a:f_bestimmt_g})}{\Longrightarrow} \, g(y^\dagger) \wedge x^\dagger = (0,0)
\,\Longleftrightarrow \,g^{\dagger}(y) \wedge x =(0,0)
\end{equation*}
and
\begin{equation*}
  g^{\dagger}(y) \wedge x = (0,0)
\, \Longleftrightarrow \, g(y^\dagger)^\dagger \wedge x = (0,0)
\, \stackrel{(\ref{eq:a:g_bestimmt_f})}{\Longrightarrow}  \,  y \wedge f(x) = (0,0)
\end{equation*}
``$\Leftarrow$'': Assume the right hand side of (\ref{dbm}) holds. Then for $x,y \in \{0,1\}^2$
\begin{equation*}
y \wedge f(x)^\dagger = (0,0) \, \Leftrightarrow \, y^\dagger \wedge f(x) =(0,0)\, \Leftrightarrow \, g^{\dagger}(y^\dagger) \wedge x = (0,0) \, \Leftrightarrow \, g(y) \wedge x^\dagger = (0,0)
\end{equation*}
which means that (\ref{eq:a:f_bestimmt_g}) and (\ref{eq:a:g_bestimmt_f}) hold showing that $f$ and $g$ are dual.
\end{proof}
We now collect some consequences arising from this characterisation of duality.
\begin{lemma}\label{eigdbm}
Let $f,g$ be basic mechanisms, which are dual. Then the following statements hold:
\begin{itemize}
\item[(i)] $f(0,0) = (0,0)$ and thus $g(0,0)=(0,0)$.
\item[(ii)] $f$ and $g$ are monotonous, i.e. it holds that
$$(\forall x,y \in \{0,1\}^2) \, x \leq y \, \Longrightarrow \, f(x)\leq f(y) \text{ and } g(x)\leq g(y)$$ 
\end{itemize}
The relation ``$\leq$'' is interpreted component-wise.
\end{lemma}
\begin{proof}
Let $f,g$ be dual basic mechanisms.\\
\textit{Ad (i):} It holds by definition
$$(\forall y \in \{0,1\}^2) \, (0,0) \wedge g^{\dagger}(y) = (0,0)$$
Therefore, our characterisation (\ref{dbm}) yields 
$$(\forall y \in \{0,1\}^2) \, f(0,0) \wedge y = (0,0),$$
and thus $f(0,0)=(0,0)$.\\
\textit{Ad (ii):} As one can easily check, the following equivalence holds for $x,y \in \{0,1\}^2$.
\begin{equation}
x \leq y \, \Longleftrightarrow \, \rub{(\forall z \in \{0,1\}^2)\, y \wedge z = (0,0) \, \Rightarrow \, x \wedge z = (0,0)}
\end{equation}
Let now $x \leq y$, then by (\ref{dbm}) it holds for arbitrary $z \in \{0,1\}^2$
$$f(y) \wedge z = (0,0) \, \Rightarrow \, g^{\dagger}(z) \wedge y =(0,0) \, \stackrel{x\leq y}{\Longrightarrow} \, g^{\dagger}(z) \wedge x = (0,0) \, \Rightarrow \, f(x)\wedge z = (0,0).$$
Thus, $f(x) \leq f(y)$.
\end{proof}

One can check by direct computation that each pair $(f,g)$ of basic
mechanisms given in each row of
Figure~\ref{f:dual_mechanisms} is dual in the sense of (\ref{dbm}).
\begin{figure}[ht]
\begin{center}
\begin{tabular}{l||c|c|c|c||c|c|c|c|}
$N^o$ & $f(0,0)$ & $f(0,1)$ & $f(1,0)$ & $f(1,1)$ & $g(0,0)$ & $g(0,1)$ & $g(1,0)$ & $g(1,1)$ \\ \hline \hline
i) & (0,0) & (0,0) & (1,1) & (1,1) & (0,0) & (0,1) & (0,1) & (0,1) \\ \hline
ii) & (0,0) & (0,1) & (1,1) & (1,1) & (0,0) & (0,1) & (1,1) & (1,1) \\ \hline
iii) & (0,0) & (0,0) & (0,1) & (0,1) & (0,0) & (0,0) & (0,1) & (0,1) \\ \hline
iv) & (0,0) & (0,1) & (1,0) & (1,1) & (0,0) & (0,1) & (1,0) & (1,1) \\ \hline
v) & (0,0) & (0,0) & (0,0) & (0,0) & (0,0) & (0,0) & (0,0) & (0,0) \\ \hline
vi) & (0,0) & (1,1) & (1,1) & (1,1) & (0,0) & (1,1) & (1,1) & (1,1) \\ \hline
\end{tabular}
\end{center}
\caption{\footnotesize Six pairs of dual mechanisms}
\label{f:dual_mechanisms}%
\end{figure}
From those six dualities, other dualities can be derived by using the following lemma.
\begin{lemma}   \label{l:more_dual_pairs}
Let $f,g$ be basic mechanisms. Then it holds:
\begin{itemize}
\item[(i)] $f$ and $g$ are dual iff $f^{\dagger}$ and $g^{\dagger}$ are dual.
\item[(ii)] $f$ and $g$ are dual iff $\hat{f}$ and $\hat{g}^{\dagger}$ are dual.
\end{itemize}
\end{lemma}
\begin{proof}
Let $f,g$ be basic mechanisms.\\
Using (\ref{dbm}) and the fact that $(f^{\dagger})^{\dagger} =f$ we have
\begin{equation}\begin{split}
f &\text{ and } g \text{ are dual.}\\
& \Longleftrightarrow \,  \rub{(\forall x,y \in \{0,1\})\; y \wedge f(x)= (0,0)  \, \Leftrightarrow \,  g^{\dagger}(y) \wedge x = (0,0)}\\
& \Longleftrightarrow \,  \rub{(\forall x,y \in \{0,1\})\; y^\dagger \wedge f(x)^\dagger= (0,0)  \, \Leftrightarrow \,  g^{\dagger}(y)^\dagger \wedge x^\dagger = (0,0)}\\
& \Longleftrightarrow \,  \rub{(\forall x,y \in \{0,1\})\; y \wedge f(x^\dagger)^\dagger= (0,0)  \, \Leftrightarrow \,  g^{\dagger}(y^\dagger) \wedge x = (0,0)}\\
& \Longleftrightarrow \,  \rub{(\forall x,y \in \{0,1\})\; y \wedge f^{\dagger}(x)= (0,0)  \, \Leftrightarrow \,  g(y) \wedge x = (0,0)}\\
&  \Longleftrightarrow \, f^{\dagger} \text{ and } g^{\dagger}  \text{ are dual.}
\end{split}
\end{equation}
This proves assertion (i). Considering the second assertion, we obtain by using (\ref{dbm}) 
\begin{equation}\begin{split}
\hat{f} &\text{ and } \hat{g}^{\dagger} \text{ are dual.}\\
& \Longleftrightarrow \,  \rub{(\forall x,y \in \{0,1\})\; y \wedge \hat{f}(x)= (0,0)  \, \Leftrightarrow \,  \hat{g}(y) \wedge x = (0,0)}\\
& \Longleftrightarrow \,  \rub{(\forall x,y \in \{0,1\})\; y \wedge f(x^\dagger)= (0,0)  \, \Leftrightarrow \,  g(y^\dagger) \wedge x = (0,0)}\\
& \Longleftrightarrow \,  \rub{(\forall x,y \in \{0,1\})\; y \wedge f(x^\dagger)= (0,0)  \, \Leftrightarrow \,  g^{\dagger}(y) \wedge x^\dagger = (0,0)}\\
& \Longleftrightarrow \,  \rub{(\forall x,y \in \{0,1\})\; y \wedge f(x)= (0,0)  \, \Leftrightarrow \,  g^{\dagger}(y) \wedge x = (0,0)}\\
&  \Longleftrightarrow \, f \text{ and } g \text{ are dual.}
\end{split}
\end{equation}
which proves assertion (ii).
\end{proof}

It remains to show that  all pairs of dual basic mechanisms are given by
Figure~\ref{f:dual_mechanisms} together with Lemma~\ref{l:more_dual_pairs}.
The number of maps $f:\{0,1\}^2 \rightarrow \{0,1\}^2$ is $4^4=256$.
By Lemma \ref{eigdbm}, a basic mechanism $f$ which has a dual must satisfy
$f(0,0)=(0,0)$.
This observation reduces the number of possible basic mechanisms with a dual to $4^3=64$.
Taking into account the monotonicity from Lemma \ref{eigdbm},
the basic mechanisms $f$ which have to be investigated further
are listed in Figure~\ref{f:remaining_mechanisms}.
There are $25=1+4+4+16$ basic mechanisms left;
one with $f(1,1)=(0,0)$,
four with $f(1,1)=(1,0)$,
four with $f(1,1)=(0,1)$
and all $4\times 4$  with $f(1,1)=(1,1)$.
  \newlength{\lenneu}%
  \setlength{\lenneu}{\linewidth*\real{0.9}}%
  \newlength{\drittelschub}%
  \setlength{\drittelschub}{(\linewidth-\lenneu)/3}%
  \newlength{\sechstelschub}%
  \setlength{\sechstelschub}{\drittelschub/2}%
  \newlength{\lenneuhalb}%
  \setlength{\lenneuhalb}{(\lenneu)/2}%
  \begin{figure}[ht]
  \begin{beidschub}[\drittelschub]
  \noindent
  \begin{minipage}[t]{\lenneuhalb}
     \begin{tabular}[t]{r||c|c|c|c|}
     $N^o$ & $f(0,0)$ & $f(0,1)$ & $f(1,0)$ & $f(1,1)$ \\ \hline \hline
     1 & (0,0) & (0,0) & (0,0) & (0,0) \\ \hline
     2 & (0,0) & (0,0) & (0,0) & (0,1) \\ 
     3 & (0,0) & (0,0) & (0,0) & (1,0) \\ 
     4 & (0,0) & (0,0) & (0,0) & (1,1) \\ \hline
     5 & (0,0) & (0,0) & (0,1) & (0,1) \\
     6 & (0,0) & (1,0) & (0,0) & (1,0) \\
     7 & (0,0) & (0,0) & (1,0) & (1,0) \\
     8 & (0,0) & (0,1) & (0,0) & (0,1) \\ \hline
     9 & (0,0) & (0,0) & (0,1) & (1,1) \\
     10 & (0,0) & (0,0) & (1,0) & (1,1) \\ 
     11 & (0,0) & (0,1) & (0,0) & (1,1) \\
     12 & (0,0) & (1,0) & (0,0) & (1,1) \\ \hline
     \end{tabular}
  \end{minipage}
  \hspace{\sechstelschub}
  \begin{minipage}[t]{\lenneuhalb}
     \begin{tabular}[t]{r||c|c|c|c|}
     $N^o$ & $f(0,0)$ & $f(0,1)$ & $f(1,0)$ & $f(1,1)$ \\ \hline \hline
     13 & (0,0) & (0,1) & (0,1) & (0,1) \\
     14 & (0,0) & (1,0) & (1,0) & (1,0) \\ 
     15 & (0,0) & (0,0) & (1,1) & (1,1) \\
     16 & (0,0) & (1,1) & (0,0) & (1,1) \\ \hline
     17 & (0,0) & (0,1) & (1,0) & (1,1) \\
     18 & (0,0) & (1,0) & (0,1) & (1,1) \\ \hline
     19 & (0,0) & (0,1) & (0,1) & (1,1) \\
     20 & (0,0) & (1,0) & (1,0) & (1,1) \\ \hline 
     21 & (0,0) & (0,1) & (1,1) & (1,1) \\
     22 & (0,0) & (1,1) & (1,0) & (1,1) \\
     23 & (0,0) & (1,0) & (1,1) & (1,1) \\
     24 & (0,0) & (1,1) & (0,1) & (1,1) \\ \hline
     25 & (0,0) & (1,1) & (1,1) & (1,1) \\ \hline
     \end{tabular}
  \end  {minipage}
  \end{beidschub} 
  \caption{\footnotesize List of basic mechanisms which might
  have a dual basic mechanism}\label{f:remaining_mechanisms}
  \end{figure}

\noindent We now consider those 25
 mechanisms successively:\\

\noindent\textit{Mechanism 1:}\\
Mechanism 1 is self-dual (cf. v) in the table of dualities above).\\

\noindent\textit{Mechanisms 2-4:}\\
Assume one of those mechanisms has a dual $g$. Then $f(0,1)=(0,0)$ implies using (\ref{dbm})
\begin{equation}
\begin{split}(\forall y &\in \{0,1\}^2) \, y \wedge f(0,1) = (0,0)\\ & \Longrightarrow \,(\forall y \in \{0,1\}^2)\, g(y^\dagger) \wedge (0,1) = (0,0) \\
& \Longrightarrow\, (\forall x \in \{0,1\}^2) \, g(x) \in \{(0,0),(1,0)\}.
\end{split}
\end{equation}
In the same way $f(1,0)=(0,0)$ implies $ (\forall x \in \{0,1\}^2) \, g(x) \in \{(0,0),(0,1)\}$. This yields that $g \equiv (0,0)$. Therefore, $g(1,1) \wedge (1,1) = (0,0)$ but $(1,1) \wedge f(1,1) \neq (0,0)$, which means that (\ref{dbm}) does not hold in contradiction to the assumption that $g$ is a dual.\\
So the mechanisms $2-4$ do not have a dual.\\

\noindent\textit{Mechanisms 5-8:}\\
Mechanism 5 is self-dual, as given under iii) in the table of dualities. Denote mechanism 5 by $f$. Then Mechanism 6 is $f^{\dagger}$ and therefore also self-dual. Mechanism 8 is equal to $\hat{f}$
and mechanism 7 is equal to $\hat{f}^{\dagger}$.
Therefore, the mechanisms 7 and 8 are dual to each other.

\bigskip
\noindent\textit{Mechanisms 9-12:}\\
Assume that mechanism 9 has a dual $g$. $f(0,1)= (0,0)$ then implies that $ (\forall x \in \{0,1\}^2) \, g(x) \in \{(0,0),(0,1)\}$. On the other hand $f(1,0) = (0,1)$ yields
\begin{equation}
\begin{split}(\forall y &\in \{(0,0),(1,0)\}) \, y \wedge f(1,0) = (0,0)\\ & \Longrightarrow \,(\forall y \in \{(0,0),(1,0)\})\, g(y^\dagger) \wedge (0,1) = (0,0) \\
& \Longrightarrow\, g(0,1) \in \{(0,0),(1,0)\}.
\end{split}
\end{equation}
Thus, $g(0,1)=(0,0)$, which leads to $g(0,1) \wedge (1,1) = (0,0)$ while $(1,0) \wedge f(1,1) \neq (0,0)$ in contradiction to the assumption that $g$ is a dual mechanism.\\
Denote mechanism 9 by $f$.
The mechanisms 10, 11 and 12 are equal to $\hat{f}^{\dagger}$, $\hat{f}$
and $f^{\dagger}$, respectively. Thus, neither of the mechanisms 10-12 has
a dual mechanism.

\bigskip
\noindent\textit{Mechanisms 13-16:}\\
The mechanisms 13 and 15 are dual to each other (cf. i) in the table of dualities). Mechanism 14 is the transpose of mechanism 13 and mechanism 16 is the transpose of mechanism 15. Therefore, the mechanisms 14 and 16 are dual to each other.\\

\noindent\textit{Mechanisms 17-18:}\\
The mechanism 17 is self-dual (cf.  iv) in the table of dualities).
Denote mechanism 17 by $f$. Mechanism 18 is both equal to $\hat{f}$ 
and to $\hat{f}^{\dagger}$ are therefore self-dual, too.

\bigskip
\noindent\textit{Mechanisms 19-20:}\\
Assume that mechanism 19, denoted by $f$, has a dual $g$. Then as above $f(1,0)=(0,1)$ implies that $ g(0,1) \in \{(0,0),(1,0)\}$. On the other hand $f(0,1)=(0,1)$ yields
\begin{equation}
\begin{split}(\forall y &\in \{(0,0),(1,0)\}) \, y \wedge f(0,1) = (0,0)\\ & \Longrightarrow \,(\forall y \in \{(0,0),(1,0)\})\, g(y^\dagger) \wedge (1,0) = (0,0) \\
& \Longrightarrow\, g(0,1) \in \{(0,0),(0,1)\}.
\end{split}
\end{equation}
This means $g(0,1)=(0,0)$ which together with the fact that $f(1,1)=(1,1)$
again leads to a contradiction.\\
Mechanism 20 is $f^{\dagger}$. Thus, mechanism 20 has no
dual, too.

\bigskip
\noindent\textit{Mechanisms 21-24:}\\
Mechanism 21, denoted by $f$, is self-dual (cf. ii) in the table of dualities) and mechanism 22 is $f^\dagger$ and thus also self-dual.
The mechanisms 23 and 24 are equal
to $\hat{f}^{\dagger}$ and $\hat{f}$, respectively, and therefore dual
to each other.\\

\noindent\textit{Mechanism 25:}\\
This mechanism is self-dual (cf. vi) in the table of dualities).\\

\bigskip\noindent
Thus, there are $16$ basic mechanisms which have a dual, $8$ of which are self-dual.


\hyphenation{Sprin-ger}

\end{document}